\documentclass{article}

\usepackage[latin1]{inputenc}
\usepackage{amsmath,amssymb,amsthm}
\usepackage{graphicx,graphics,color,epsfig}
\usepackage{hyperref}

\newcommand{\veps}{\varepsilon}

\newcommand{\mpdiv}{\text{$/$\hspace{-1.42mm}\tiny${}^\circ $\normalsize}}

\newcommand{\CC}{\mathcal C}

\newcommand{\CG}{\mathcal G}

\newcommand{\CP}{\mathcal P}

\newcommand{\CS}{\mathcal S}

\newcommand{\CU}{\mathcal U}
\newcommand{\CV}{\mathcal V}
\newcommand{\CW}{\mathcal W}
\newcommand{\CX}{\mathcal X}

\newcommand{\BE}{\mathbb E}

\newcommand{\BN}{\mathbb N}

\newcommand{\BR}{\mathbb R}

\theoremstyle{plain}
\newtheorem{theorem}{Theorem}
\newtheorem{proposition}{Proposition}
\newtheorem{corollary}{Corollary}

\theoremstyle{definition}

\theoremstyle{remark}

\title{\textsc{Fundamental Diagrams of 1D-Traffic Flow by Optimal Control
        Models}}
\author{Nadir Farhi\\ \vspace{-2mm}
\small{INRIA - Paris - Rocquencourt}\\ 
\small{Domaine de Voluceau, 78153, Le Chesnay, Cedex France.}\\
\small{\texttt{nadir.farhi@inria.fr}}}

\date{}

\begin{document}

\maketitle

\begin{abstract}
Traffic on a circular road is described by dynamic programming
equations associated to optimal control problems. By solving the
equations analytically, we derive the relation between the average
car density and the average car flow, known as the
\emph{fundamental diagram of traffic}. First, we present a
model based on min-plus algebra, then
we extend it to a stochastic dynamic programming model, then
to a stochastic game model. The average car flow is derived as the
average cost per time unit of optimal control problems, obtained
in terms of the average car density. The models presented in this
article can also be seen as developed versions of the car-following
model. The derivations proposed here can be used to approximate, understand and
interprete fundamental diagrams derived from real measurements.
\end{abstract}


\section{Introduction}

The relation between the car density and the car flow in road traffic systems is known under the name of
\emph{fundamental diagram of traffic}. Basically, fundamental diagrams are studied on one road (a urban road, a highway segment
or a circular ring)~\cite{Green35, Der94, FY96-1, Bla00, WWHH00, Hel01, Dag05, LMQ05, Far08}. In this case of one road without crossing, one talks about fundamental diagram of one dimensional traffic (1D-traffic). However, many works on fundamental diagrams of 2D-traffic (roads with crossings) have appeared recently~\cite{NS92, FY96-2, CSS00, FY01-1, FY01-2, GD07, DG08, GD08, Far08, BL09, Hel09}. We are intereseted in this article by fundamental diagrams of 1D-traffic. We present an optimal control approach that permits to understand, approximate and interpret 1D-diagrams.

The relation flow-density have been observed on a highway since 1935 by Greenshields~\cite{Green35}. Lighthill, Whitham, and Richards (LWR)~\cite{LW55, Ric56} describe the traffic by a car conservation equation $\partial_t \sigma+\partial_x \rho=0$, where $\sigma(x,t)$ and $\rho(x,t)$ denote respectively the density and the flow of vehicles in position $x$ at time $t$. In the stationary regime, the flow $\rho$ is linked to the density $\sigma$ by the following functional relation: $\rho(\sigma)=d \bar{v}$, where $\bar{v}$ is the average car-speed. To complete the dynamics given by the car conservation equation, LWR supposed the existence of a traffic behavior equation $\rho=f(\sigma)$, either in the situation of time and space dependence, called the fundamental diagram of traffic.

A well-known microscopic model is the car-following model~\cite{HMPR59,GHR61}. The traffic is described on one road where it is supposed that vehicles follow their predecessors without overtaking and with a stimulation-response relation. In~\cite{DG08}, Daganzo and Gerolimins used a variational theory~\cite{Dag05, GD07, GD08} to show the existence of a concave
macroscopic fundamental diagram on a ring~\footnote{This approach has been extended, in the same article~\cite{DG08}, to
a network, by using an aggregation method.}. Our results are very close to those given in~\cite{DG08}. However,
the approaches, the models, and the exploitation of the results are very different.

In this article, the traffic on a circular road (ring) is described by dynamic programming
equations (DPE) of optimal control problems. The average car flow is derived as the average cost
per time unit of the optimal control problem considered. In addition, the average car flow
is given in term of the average car density, giving thus the fundamental diagram of traffic.
The models we present here lead to piecewise linear diagrams.

We consider $n$ vehicles moving on a one-lane circular road without overtaking.
We start with a very simple model (the min-plus model) witch we extend by refining the traffic
description. The min-plus~\footnote{A short review in min-plus algebra is given in
section~\ref{sec-mod1}.} linear model describes the traffic using two parameters: a desired
velocity $v$, fixed and common for all vehicles, and a safety distance $\sigma$ between two successive
vehicles. The dynamics tells simply that at each time, each car tries to move with a velocity $v$
under the constraint that it has to leave a safety distance $\sigma$ with respect to the car ahead. The
dynamics are given by a min-plus liear system~\footnote{which can be seen as a dynamic programming equation 
of a deterministic optimal control problem.}, and the average car flow is derived as the min-plus
eigenvalue of this system.

We extend the min-plus model by assuming that the desired car velocities are not constant but depend 
on the distances between successive cars. The first extension gives a model that describes the traffic
by a DPE of a stochastic optimal control problem. We solve analytically this
equation and get the fundamental traffic diagram. This extension permits to realize a large
class of cancave fundamental diagrams. The second extension gives a model that describes the traffic by a
DPE of a stochastic game problem with two players. Similarly, we solve
analytically the DPE and get the fundamental traffic diagram. This latter extension permits to realize
even non concave diagrams.

\section{Min-plus Traffic Model}
\label{sec-mod1}

We present in this section the first model which we call the min-plus traffic model. It is a very basic
model, presented mainly to introduce its extensions. This model is a dual version of the min-plus
model studies in~\cite{LMQ05}.
Let us first give a short review of the min-plus
algebra. \emph{Min-plus algebra}~\cite{BCOQ92} is the commutative idempotent semiring
$(\BR\cup\{+\infty\},\oplus,$ $\otimes)$ where the operations $\oplus$ and $\otimes$ are defined by $a\oplus
b=\min(a,b)$ and $a\otimes b=a+b$ respectively. We denote this structure by $\BR_{\min}$. The zero and the
unity elements are repectively $+\infty$ denoted $\veps$ and $0$ denoted $e$. The main differences between
standard and min-plus algebras are the idempotency ($a\oplus a=a,\;\forall a\in\BR_{\min}$) and the non
simplification ($a\oplus b=a\oplus c \nRightarrow b=c$) of min-plus addition. The structure $\BR_{\min}$
on scalars induces another idempotent semiring on the set of square matrices with entries in $\BR_{\min}$.
If $A$ and $B$ are two square matrices with entries in $\BR_{\min}$ (we say min-plus square matrices), then
the addition is defined by: $(A\oplus B)_{ij}=A_{ij}\oplus B_{ij}, \; \forall i,j$, and the product by~:
$(A\otimes B)_{ij}=\bigoplus_{k}(A_{ik}\otimes B_{kj}), \; \forall i,j$. The zero and the unity matrices
are also denoted by $\veps$ and $e$ respectively.
A directed graph $\CG(A)$ is associated to a square min-plus matrix $A$. It is the graph whose nodes correspond
to the matrix lines and whose arcs correspond to the no null ($\neq \veps$) entries of $A$. When
$A_{ij}\neq \veps$, there exists an arc in $\CG(A)$ going from node $j$ to node $i$.
\begin{theorem}\cite{BCOQ92}\label{spec}
  If the graph $\CG(A)$ associated to a min-plus square matrix $A$ is strongly connected, then
  $A$ admits a unique min-plus eigenvalue $\mu$ given by the minimum of the average weights of the graph circuits:
  $\mu=\min_{c\in \CC}(|c|_w/|c|_l)$,
  where $\CC$ is the set of the circuits in $\CG(A)$, $|c|_w$ is the weight of a circuit $c$ given by the
  min-plus product (standard sum) of the arc weights, and $|c|_l$ is the circuit length given by the number
  of arcs of the circuit.
\end{theorem}
\begin{theorem}\cite{BCOQ92}\label{period}
  The min-plus linear dynamic system associated to a square min-plus matrix $A$ whose graph
  is strongly connected, defined by:
  $x^{k+1}=A\otimes x^k$,
  is asymptotically periodic:
  $\exists T,K,\mu: \forall k\geq K:A^{k+T}=\mu^T\otimes A^k\;$.
  Moreover, $\mu$ coincides with the unique eigenvalue of~$A$.
\end{theorem}

\subsection{The model}

\begin{figure}[h]
  \begin{center}
    \includegraphics[width=3cm,height=2.7cm]{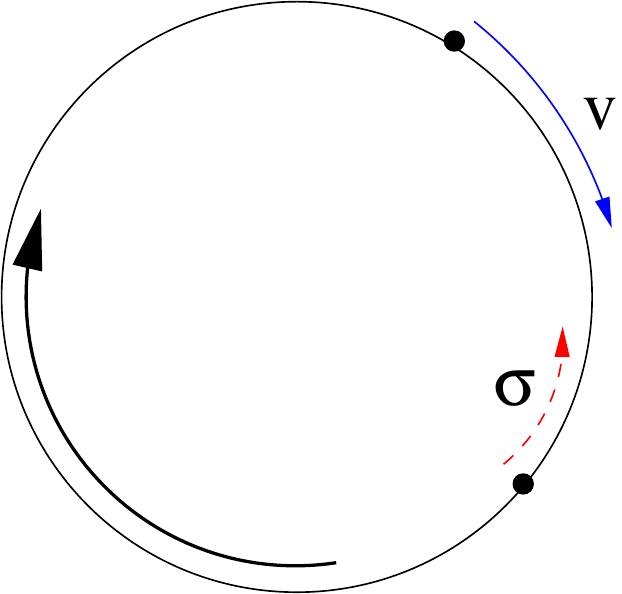}
    \caption{Traffic on a circular road.}
    \label{circul}
  \end{center}
\end{figure}

We assume here that all cars have one length, and we take this length as the unity of distance. We consider
$n$ car moving on a one-lane circular road of length $m$ (the road cannot contain more than $m$ cars), with
$n\leq m$; see Figure~\ref{circul}. The car density on the road is $n/m$. We assume that the cars have
one same desired velocity $v$ and that each car has to respect a safety distance $\sigma$ with respect to
the car ahead.

\begin{figure}[h]
  \begin{center}
    \includegraphics[width=4.5cm,height=4.5cm]{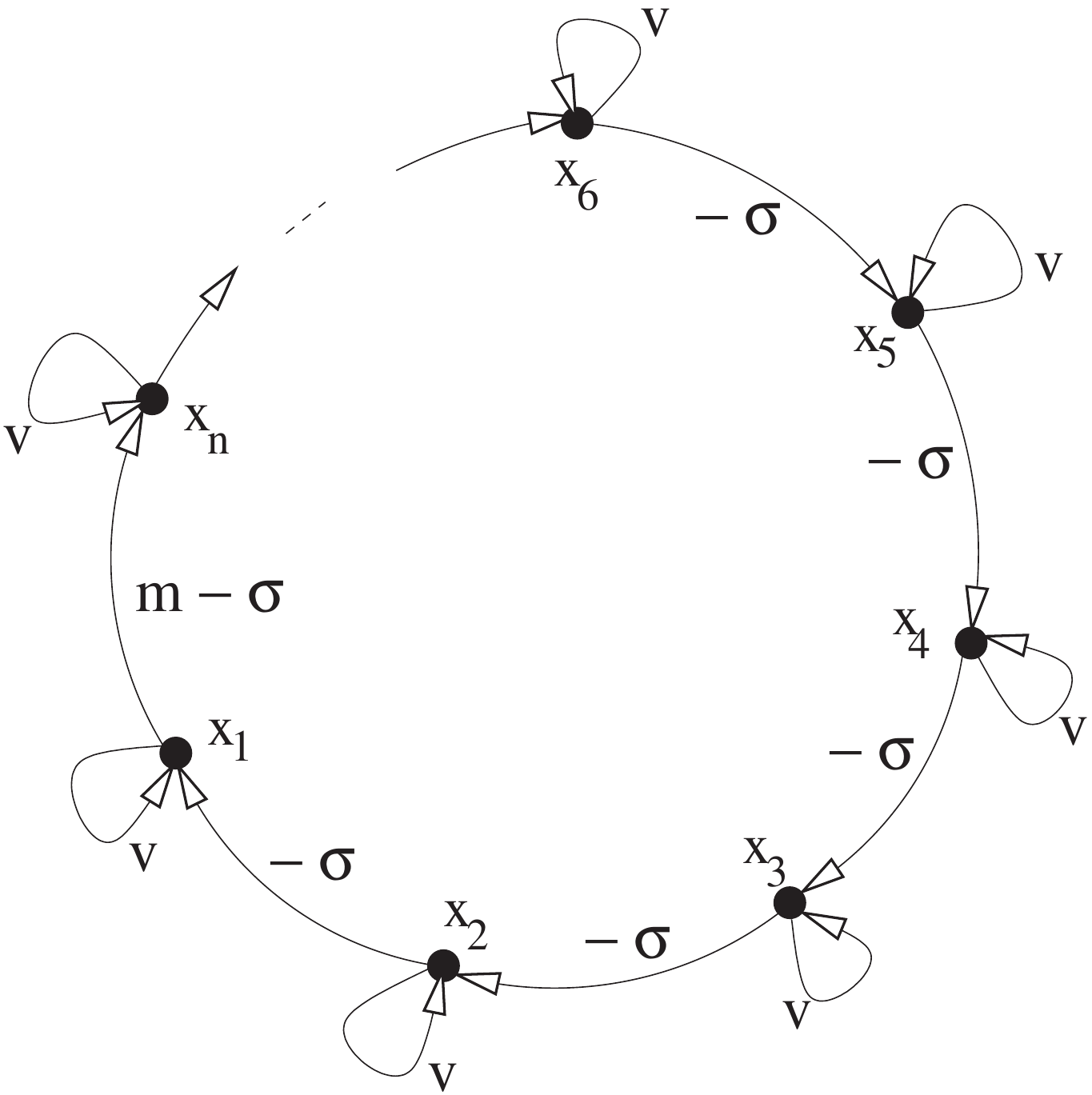}
    \caption{The graph associated to $A$.}
    \label{graph4}
  \end{center}
\end{figure}

Let us denote by $x_i^k$ the distance travelled by a car $i$ up to time $k$. This distance satisfies the following dynamics:
\begin{equation}\label{dis1}
  x_i^{k+1}=\begin{cases}
              \min\{v+x_i^k, x_{i+1}^k-\sigma\} \quad \text{ if } i<n,\\
              \min\{v+x_i^k, x_1^k+m-\sigma\} \quad \text{ if } i=n.
            \end{cases}
\end{equation}
The average growth rate per time unit of system (\ref{dis1}) is interpreted as the average car velocity on the road.
This system is written in min-plus algebra as follows~:
\begin{equation}\label{dis2}
  x_i^{k+1}=\begin{cases}
        v \otimes x_i^k\oplus (e\mpdiv\sigma) \otimes x_{i+1}^k \quad \text{ if } i<n,\\
        v \otimes x_i^k\oplus (m\mpdiv\sigma) \otimes x_1^k \quad \text{ if }  i=n,
            \end{cases}
\end{equation}
where the symbol $\mpdiv$ denotes the standard substraction. For example, $e\mpdiv \sigma$ is nothing but $(-\sigma)$ in
standard algebra. The dynamics (\ref{dis2}) is  min-plus linear and can be written~:
\begin{equation}\label{sys1}
  x^{k+1}=A\otimes x^k,
\end{equation}
where $A$ is a min-plus matrix given by~:
$$A=\begin{bmatrix}
       v & e\mpdiv\sigma &\veps & \cdots & \veps\\
       \veps & v & e\mpdiv\sigma & \cdots & \veps\\
       \vdots &  & \ddots & \ddots & \vdots\\
       \veps & \veps & \veps & & e\mpdiv\sigma\\
       m\mpdiv\sigma & \veps & \veps & & v
    \end{bmatrix}.$$
\begin{theorem}\label{th-rapp}\cite{LMQ05}
  There exists an average car velocity $\bar{v}$. It is the eigenvalue of the matrix $A$ associated to
  system~(\ref{sys1}), and is given by~: $\bar{v}=\min\{v,(m-n\sigma)/n\}$.
\end{theorem}
\proof The graph associated to the min-plus matrix $A$ is shown on Figure~\ref{graph4}.
Theorem~\ref{spec} gives the eigenvalue of $A$ as the minimum of the average weights of the circuits of the graph
associated to $A$. The elementary circuits of the graph of Figure \ref{graph4} are~:
\begin{itemize}
  \item the circuit passing by all the nodes, with an average weight of $(m-n\sigma)/n$,
  \item the loops of weight $v$.
\end{itemize}
Thus we obtain $\bar{v}$. The average growth rate per time unit of system~(\ref{sys1}) is interpreted as the average
car speed on the road. Using Theorem \ref{period}, we conclude that the average car speed is the eigenvalue $\bar{v}$ of~$A$. \endproof
\begin{corollary}\label{cor1}\cite{LMQ05}
  The fundamental traffic diagram on a circular road where the traffic is described by the dynamics~(\ref{dis1}) is~:
  $f=\min\{vd,1-\sigma d\}$.
\end{corollary}
\proof
  We know that the average car flow $f$ is given by the average car speed $\bar{v}$ multiplied by the car density $d$~:   
  $f=d\bar{v}$. By replacing $\bar{v}$ by its value given in Theorem~\ref{th-rapp}, we obtain the result. \endproof

\section{Stochastic optimal control model}
\label{socm}

In this section, we extend the min-plus traffic model given in the preceding section. We assume here that each car
chooses its desired velocity depending on the distance with respect to the car ahead. The car dynamics will be written as a
dynamic programming equation of a stochastic optimal control problem. The fundamental traffic diagram is then 
derived by solving this equation. The extension we give in this section permits to approximate a large class of 
emph{concave} fundamental diagrams.

In order to clarify the modeling context and to set notations, let us give a short review on stochastic optimal control
of Markov chains. A stochastic optimal control problem in finite horizon with undiscounted costs  is written as follows
(see for example~\cite{Whi86})~:
\begin{equation}\label{socp}
  \min_{s\in\CS}\BE\left\{\lim_{T\to +\infty}\frac{1}{T}\sum_{k=0}^{T-1}c^{u^k}_{x^k}\right\}
\end{equation}
where $(x^k)_{k\in\BN}$ is a controlled Markov chain with a finite
set of states $\CX=\{1,2,\cdots,n\}$, $u^k\in\CU$ is the decision
variable taken at time $k$, with $\CU$ a finite set of controls,
$c^{u^k}_{x^k}\in\BR$ is the cost to pay at time $k$ being in $x^k$
and taking the decision $u^k$, and $\CS$ is the set of control strategies, that is the set of time-indexed sequences
in $\CU$ (a strategy $s\in\CS$ is a fixed sequence $\{u^k\}_{k\in\BN}, \; u^k\in\CU$)~\footnote{It is known that
solving optimization problem~(\ref{socp}) in $\CS$ is equivalent to solve the same problem in $\CP\subset \CS$, where
$\CP$ is the set of feedback strategies defined on $\CX$. A feedback $p$ associates to each state $x\in\CX$ a 
control $u\in\CU$ ($\CP\ni p:\CX\ni x\mapsto u\in\CU$).}.

If we denote by $M^u, u\in\CU$, the transition matrix of the Markov
chain associated to a control $u\in\CU$, then the stochastic
dynamic programming equation associated to problem (\ref{socp})
is written~:
\begin{equation}\label{epd-ergo-sto}
  \mu+v_x=\min_{u\in\CU}\{[M^u v]_x+c^u_x\}, \quad \forall\; 1\leq x\leq n,
\end{equation}
In (\ref{epd-ergo-sto}), $\mu$ can be seen as an additive eigenvalue
associated to the eigenvector $v$ of an operator $h: \BR^n\ni
v\mapsto h(v)\in\BR^n$ defined by:
\begin{equation}\label{oph}
   h_x(v)=\min_{u\in\CU}\{[M^u v]_x+c^u_x\}, \quad \forall\; 1\leq x\leq n.
\end{equation}
Operator $h$ given in~(\ref{oph}) is additive 1-homogeneous (that is
$h(\mu+v)=\mu+h(v),\;\forall\; \mu\in\BR,\; \forall\; v\in\BR^n$),
monotone ($\forall v, w\in \BR^n,\; \left[v_x\leq w_x,\; \forall\;
x\right]\Rightarrow \left[h_x(v)\leq h_x(w),\; \forall x\right]$)
and concave ($\forall v,w\in\BR^n,\;\forall t\in[0,1],
h(tv+(1-t)w)\geq th(v)+(1-t)h(w)$).

Let us define, as in~\cite{GG99}, an oriented graph $\CG(h)$ associated to an additive
1-homogeneous and monotone map $h$ by the graph of $n$ nodes where
arcs are determined as follows: there exists an arc from $i$ to
$j$ if $lim_{\nu\to\infty}f_j(\nu e_i)=\infty$, where $e_i$
denotes the $i^{th}$ vector of the canonic basis of $\BR^n$. We
denote by $\chi(h)$ the average growth rate per time unit of
the dynamical system: $v^{k+1}=h(v^k)$, defined by:
$\chi(h)=\lim_{k\to\infty}v^k/k$. In the following, we recall an important
result on additive 1-homogeneous and monotone maps.
\begin{theorem}\cite{GK95,GG98a}
  If $h$ is an additive 1-homogeneous and monotone map, and if
  $\CG(h)$ is strongly connected, then the additive eigenvalue problem
  $\mu+v=h(v)$ admits a solution $(\mu,v)$, where $v$ is
  defined up to an additive constant, not necessarily in a unique
  way, and $\mu$ is unique and satisfies $\chi(h)={}^t(\mu,\mu,\cdots,\mu)$.
\end{theorem}
\begin{corollary}\label{coro}
  Let $h$ be the operator defined in (\ref{oph}). If  $\CG(h)$ is strongly
  connected, then the eigenvalue problem (\ref{epd-ergo-sto}) admits a solution $(\mu,v)$
  where $v$ is defined up to an additive constant, not necessarily in a unique
  way, and $\mu$ is unique and satisfies $\chi(h)={}^t(\mu,\mu,\cdots,\mu)$.
\end{corollary}

\subsection{The model}

As above, we suppose $n$ vehicles moving on a circular road of length $m$, with $n\leq m$.
Let us denote by $y^k$ the distance travelled by a given vehicle
up to time $k$ and by $z^k$ the distance travelled by the vehicle
ahead up to time $k$. We add to the two constraints of velocity
limitation and safety distance, another constraint which expresses
the dependence of the velocity at time $k$ on the distance
$z^k-y^k$. Thus we obtain three constraints~:
\begin{itemize}
  \item Velocity limitation:
    $$y^{k+1}\leq y^k+v\;.$$
  \item Safety distance:
    $$y^{k+1}\leq z^k-\sigma\;.$$
  \item Dependence of the velocity on the distance $z^k-y^k$:
    $$y^{k+1}\leq y^k+\beta(z^k-y^k),\quad 0\leq \beta \leq 1\; .$$
\end{itemize}
These three constraints can be summarized in:
\begin{equation}\label{gconst}
  y^{k+1}\leq y^k+\alpha+\beta(z^k-y^k), \quad 0\leq\beta\leq 1\;.
\end{equation}
Indeed, the first constraint is obtained by taking $\alpha=v$, and
$\beta=0$, the second one by taking $\alpha=-\sigma$, and $\beta=1$,
and the third one by taking $\alpha=0$.
In general, we assume that each vehicle has to satisfy a set $\CU$
of traffic constraints of type~(\ref{gconst}). With $n$ vehicles
indexed by $i$, moving on a road of length $m$ (the car density is $d=n/m$), we denote by
$x_i^k$ the distance travelled by a vehicle $i$ up to time $k$.
The car dynamics is then written as follows~:
\begin{equation}\label{car-f}
  x_i^{k+1}=\begin{cases}
   \min_{u\in\CU}\{x_i^k+\alpha_u+\beta_u(x_{i+1}^k-x_i^k)\} & \text{if }\; i<n,\\
   \min_{u\in\CU}\{x_n^k+\alpha_u+\beta_u(m+x_1^k-x_n^k)\} & \text{if }\; i=n,
  \end{cases}
\end{equation}
and since $m=n/d$~\footnote{The case $d=0$ is trivial since it
corresponds to $n=0$. This case is implicitly neglected here.}, we obtain~:
\begin{equation}\label{eqt1}
  x_i^{k+1}=\begin{cases}
          \min_{u\in\CU}\{\alpha_u+(1-\beta_u)x_i^k+\beta_u x_{i+1}^k\} & \text{if }\; i<n,\\
          \min_{u\in\CU}\{\alpha_u+n\beta_u/d+(1-\beta_u)x_n^k+\beta_u x_1^k\} & \text{if }\; i=n.
            \end{cases}
\end{equation}
Let us define the matrices $M^u$ and the vectors $c^u$ for $u\in\CU$ by~:
$$M^u=\begin{bmatrix}
        1-\beta_u & \beta_u & 0 & \cdots & 0\\
        0 & 1-\beta_u & \beta_u & & 0\\
        \vdots & & \ddots & \ddots & 0\\
        0 & \cdots & \cdots & 1-\beta_u & \beta_u\\
        \beta_u & 0 & \cdots & 0 & 1-\beta_u
      \end{bmatrix},$$
\vspace{5mm}
$$c^u={}^t[\alpha_u,\;\; \alpha_u,\; \cdots,\; \alpha_u,\;\; \alpha_u+n\beta_u/d].$$
Equations~(\ref{eqt1}) are then written~:
\begin{equation}\label{gen-1}
  x_i^{k+1}=\min_{u\in\CU}\{[M^ux^k]_i+c^u_i\},\quad 1\leq i\leq n\;.
\end{equation}
System~(\ref{gen-1}) is a backward
dynamic programming equation of a stochastic optimal control
problem of a Markov chain with transition matrices
$M^u, u\in\CU$ and costs $c^u, u\in\CU$.

Let us denote by $h$ the operator giving the dynamics~(\ref{gen-1}), that is $h:\BR^n\to\BR^n$ given by~:
$$h_i(x^k)=\min_{u\in\CU}\{[M^ux^k]_i+c^u_i\},\quad 1\leq i\leq n\;.$$
\begin{proposition}\label{prop}
  The graph $\CG(h)$ associated to $h$ is strongly connected if and only if there exists $u\in\CU$ such
  that $\beta_u\neq 0$ (that is $\beta_u\in (0,1]$).
\end{proposition}
\proof ~
\begin{itemize}
  \item If $\exists u\in\CU$, such that $\beta_u\in(0,1]$, then for all $1\leq i\leq n$,
    there exists an arc on $\CG(h)$ going from $i+1$ to $i$ (modulo $n$).
    Indeed, we have:
    $$x_i^{k+1}=(1-\beta_u)x_i^k+\beta_u x_{i+1}^k+\alpha_u,$$
    and since $\beta_u>0$, we get:
    $$\lim_{\nu\to\infty}h_i(\nu e_{i+1})=\lim_{\nu\to\infty}[\beta_u \nu+\alpha_u]=\infty.$$
    where $e_{i+1}$ denotes the $(i+1)^{th}$ vector of the canonic
    basis of $\BR^n$. Thus the graph $\CG(h)$ is strongly
    connected.
  \item If $\forall u\in\CU, \beta_u=0$, then we can easily check
    that all arcs of $G(h)$ are loops; so the graphe $G(h)$ is not strongly
    connected. \endproof
\end{itemize}
In the following, we suppose that there exists $u$ in $\CU$ such
that $\beta_u\in (0,1]$. In terms of traffic, this means that each
car moves by taking into account the position of the car
ahead. With this assumption, we get an additively 1-homogeneous
and monotone operator $h$, whose associated graph is strongly
connected. 

Applying Corollary~\ref{coro}, we conclude that the system~:
\begin{equation}\label{epd-app-2}
  \mu+x_i=\min_{u\in\CU}\left\{(M^u x)_i+c^u_i\right\},\quad 1\leq i\leq n\;
\end{equation}
admits a solution $(\mu,x)$ where $x$ is defined up to an
additive constant, not necessarily in a unique way, and $\mu$ is
unique and satisfies~:
\begin{equation}\label{vitmoy}
  \mu=\lim_{k\rightarrow +\infty}\frac{1}{k}\;x_i^k,\;\;\; 1\leq i \leq n\;.
\end{equation}
$\mu$ is interpreted as the average speed of cars.
\begin{theorem}\label{th-gen-2}
  System~(\ref{epd-app-2}) admits a solution $(\mu,x)$ given by:
  $$\mu=\min_{u\in\CU}\{\alpha_u+\frac{\beta_u}{d}\} \quad \text{and} \quad
  x= {}^t[0 \quad 1/d \quad 2/d \quad \cdots \quad (n-1)/d].$$
\end{theorem}
\proof It is natural to think
that the asymptotic positions $x_i$, $1\leq i\leq n$ are uniformly
distributed on the ring. This gives the eigenvector $x$. It is
also natural to think that the optimal strategy is independent on
the state $x$, because of the symmetry of the system. Let us check
this. Let $\bar{u}\in\CU$ satisfying~:
$$\mu=\min_{u\in\CU}\{\alpha_u+\frac{\beta_u}{d}\}=\alpha_{\bar{u}}+\frac{\beta_{\bar{u}}}{d}.$$
First, we can easily check that $(\mu,x)$ given by:
$$\mu=\alpha_{\bar{u}}+\frac{\beta_{\bar{u}}}{d}, \quad \text{and} \quad
  x= {}^t[0 \quad 1/d \quad 2/d \quad \cdots \quad (n-1)/d],$$
is a solution of the system~: $\mu+x=M^{\bar{u}}x+c^{\bar{u}}\;$. \\
Second, the feedback strategy $s: x_i\to \bar{u},\;\; 1\leq i \leq n$, is
optimal, because for all $i\in\{1,2,\cdots,n\}$ and for all
$u\in\CU$ we have~:
$$[M^{\bar{u}}+c^{\bar{u}}x]_i = \mu+x_i
                               \leq \alpha_u+\frac{\beta_u}{d}+x_i
                               = [M^u x+c^u]_i$$
Thus the couple $(\mu,x)$ satisfies system (\ref{epd-app-2}).
\endproof
\begin{corollary}\label{cor2}
  The fundamental diagram on the circular road where the traffic is described by dynamics~(\ref{gen-1}) is given by~:
  $f=\min_{u\in\CU}\{\alpha_u d+\beta_u\}$.
\end{corollary}
\proof The average flow is equal to the average speed given in Theorem~\ref{th-gen-2} multiplied by~$d$. \endproof

\subsection*{Remarks}

\begin{enumerate}
  \item We make here a link between the model presented above and the car-following model~\cite{HMPR59,GHR61}.
    Daganzo~\cite{Dag06} has already linked the car-following model to his variational theory based model. 
    In~\cite{HMPR59,GHR61}, the traffic is described on one road by assuming that each vehicle follows his
    predecessor without possibility of overtaking and with a stimulation-response relation. Let $x_n(t)$
    denoting the position of the $n$-th vehicle on the road, at time $t$, and $T$ denoting the reaction time
    of a driver. The acceleration $d^2x_n(t+T)/dt^2$ of the $n$-th vehicle at time $t+T$ is given by multiplying
    by $\lambda$ the response to the stimulation $dx_{n-1}(t)/dt-dx_n(t)/dt$. We write~:
    \begin{equation}\label{cfm}
      \frac{d^2x_n(t+T)}{dt^2}=\lambda\big[\frac{dx_{n-1}(t)}{dt}-\frac{dx_n(t)}{dt}\big],
    \end{equation}
    where $\lambda$ is often taken as follows:
    \begin{equation}\nonumber
      \lambda=\frac{\lambda_0[dx_n(t)/dt]^m}{[x_{n-1}(t)-x_n(t)]^l}\;,
    \end{equation}
    with $\lambda_0$ a constant and $m$ and $l$ are parameters.
    In the simple case where $m=l=0$ i.e. $\lambda=\lambda_0$, we
    obtain the linear model~:
    \begin{equation}\nonumber
      \frac{dx_n(t)}{dt}=\lambda_0[x_{n-1}(t)-x_n(t)],
    \end{equation}
    which we can write:
    \begin{equation}\label{cf-mod}
      v=\lambda_0 s+\alpha,
    \end{equation}
    where $s=x_{n-1}-x_n$ and $\alpha$ is a constant
    determined by the boundary condition $v=0$ corresponding to
    the jam state $s=s_j$.
    The velocities considered in equation~(\ref{car-f}) are nothing but what is given in~(\ref{cf-mod}).
  \item Approximating Diagrams using the formula of Corollary~\ref{cor2} is also computing Fenchel transforms (concave
    version). This is known and used in~\cite{Dag06,ABS08}. Indeed, if we denote by $\CV$ the set
    $\CV=\{\alpha_u, \; u\in\CU\}$ and define the function $g$ by:
    $$\begin{array}{llll}
        g: & \CV & \to & \BR\\
           & v=\alpha_u & \mapsto & - \beta_u\;,
      \end{array}$$
    then we obtain:
    $$f(d)=\min_{v\in\CV}\big(d v - g(v)\big) = g^*(d),$$
    where $g^*$ denotes the Fenchel transform of $g$. Thus, giving an
    approximation of a diagram is giving a finite set $\CV=\{\alpha_u,
    u\in\CU\}$ and defining the function $g$, which associates for
    each $\alpha_u, u\in\CU$ a value $\beta_u$. Graphically, this is
    giving a finite set of segments by their slopes $\alpha_u$ and
    their values at the origin $\beta_u$.

  \item Using the stochastic optimal control model given above, we obtain a large class of concave diagrams, but
    not all the concave diagrams. Indeed, every concave function $f(d)$ can be approximated, with any precision,
    by a function $h(d)=min_{u\in\CU} (\alpha_u d+\beta_u)$ with $\alpha_u\in\BR$ and $\beta_u\in\BR$, for all
    $u$ in $\CU$; but in the model, we accept only $\beta_u$ satisfying $\beta_u\in[0,1]$.

  \item The min-plus linear model is a particular case where
    $\CU=\{u_1,u_2\}$ with $(\alpha_1,\beta_1)=(v,0)$ and
    $(\alpha_2,\beta_2)=(-\sigma,1)$. In this case, the approximation is
    a piecewise linear function with two segments.
\end{enumerate}

\section{Stochastic game model}

We extend again the stochastic dynamic programming model
to obtain a stochastic game one. We assumed in the preceding
sections that in both low and high density cases, the drivers have superior bounds of speed to respect
($\leq$~inequalities), and they maximize their speed by moving with
the minimum superior bound. The extension is to suppose also the
dual situation. Indeed, in the case of high densities, the
drivers can have inferior bounds of speed to respect ($\geq$~inequalities), 
and then minimize their speed by moving with the
maximum inferior bound. This is detailed below.
With this extension, the car dynamics are interpreted in
term of stochastic games, and the fundamental traffic diagram is
obtained, as above, by solving analytically a generalized eigenvalue problem.
Moreover, even non concave diagrams can be approximated with this extension.

Let us first give a short review on stochastic games. A stochastic game problem in infinite horizon with undiscounted
costs is written:
\begin{equation}\label{socp1}
  \min\max|_{s\in\CS}\BE\left\{\lim_{T\to +\infty}\frac{1}{T}
        \sum_{k=0}^{T-1}c^{u^k w^{k}}_{x^k}\right\}
\end{equation}
where $(x^k)_{k\in\BN}$ is a controlled Markov chain with a finite
set of states $\CX=\{1,2,\cdots,n\}$, $u^k\in\CU$ is the minimizer
decision variable taken at time $k$, with $\CU$ a finite set of
controls, $w^k\in\CW$ is the maximizer decision variable taken at
time $k$, with $\CW$ a finite set of controls,
$c^{u^kw^k}_{x^k}\in\BR$ is the cost to pay at time $k$ being in
$x^k$ and when the minimizer takes the decisions $u^k$ and the
maximizer takes the decision $w^k$, and $\CS$ is the set of control strategies, that is the set of time-indexed sequences
in $\CU\times \CW$ (a strategy $s\in\CS$ is a fixed sequence $\{(u^k,w^k)\}_{k\in\BN}, \; u^k\in\CU, w^k\in\CW$).

If we denote by $M^{uw}, u\in\CU, w\in \CW$, the transition matrix
of the controlled Markov chain associated to the controls
$u\in\CU$ and $w\in\CW$, then the stochastic dynamic programming
equation associated to problem (\ref{socp1}) (where the
maximizer knows, at each step, the choice of the minimizer) is
written:
\begin{equation}\label{epd-ergo-jeu}
  \mu+v_x=\min_{u\in\CU}\max_{w\in\CW}\{[M^{uw} v]_x+c^{uw}_x\}, \quad \forall\; 1\leq x\leq n.
\end{equation}

In (\ref{epd-ergo-jeu}), $\mu$ can be seen as an additive
eigenvalue of an operator $h: \BR^n\ni v\mapsto h(v)\in\BR^n$
defined~by:
\begin{equation}\label{oph2}
  h_x(v)=\min_{u\in\CU}\max_{w\in\CW}\{[M^{uw} v]_x+c^{uw}_x\} \quad \forall\; 1\leq x\leq n.
\end{equation}
The operator $h$ is additive 1-homogeneous and monotone. Corollary~\ref{coro} can be applied again.

\subsection{The model}

Here we extend the model by taking into account the driver's behavior
changing from low densities to high ones. The
difference between these two situations is that in low densities,
drivers, moving, or being able to move with high velocities,
they try to leave large safe distances between each other, so the
safe distances are maximized; while in high densities,
drivers, moving, or having to move with low velocities, they try
to leave small safe distances between each other in order to avoid
jams, so they minimize safe distances. To illustrate this idea,
let us denote by $y^k$ (resp. $z^k$) the travelled distance up to time
$k$ by a given car (resp. by the car ahead). Instead of
maintaining the safe distance more than $\sigma$ i.e. $y^{k+1}\leq
z^k-\sigma\;$, let's use the constraint:
$$y^{k+1}\leq \max\{z^k-\sigma,(y^k+z^k)/2\}\;.$$
In a low density situation where the vehicles are separated by at
least $2\sigma$ we have~:
$$\max\{z^k-\sigma,(y^k+z^k)/2\}=z^k-\sigma\;,$$
while in a high density situation we obtain:
$$\max\{z^k-\sigma,(y^k+z^k)/2\}=(y^k+z^k)/2\;.$$
In this latter case, we accept the vehicles moving closer but by
reducing the approach speed in order to avoid collisions. The
whole dynamics will be~:
$$y^{k+1}=\min\{\max\{z^k-\sigma,(y^k+z^k)/2\},y^k+v\}\;.$$
In general, we denote by $x^k_i$ (resp. $x^k_{i+1}$) the distance travelled by
a car $i$ (resp. by the car ahead) up to time $k$. We
obtain the following dynamics~:
\begin{equation}\label{eq-jeu}
  x_i^{k+1}=\begin{cases}
    \min_{u\in\CU}\max_{w\in\CW}\{(1-\beta_{uw})x_i^k+\beta_{uw}x_{i+1}^k+\alpha_{uw}\} & \text{si } i<n\;,\\
    \min_{u\in\CU}\max_{w\in\CW}\{(1-\beta_{uw})x_n^k+\beta_{uw}x_1^k+\alpha_{uw}+n\beta_{uw}/d\} &
        \text{si } i=n\;,
  \end{cases}
\end{equation}

As in Section~\ref{socm}, we define the matrices $M^{uw},
(u,w)\in(\CU\times \CW)$ and the vectors $c^{uw},
(u,w)\in(\CU\times \CW)$~by:
$$M^{uw}=\begin{bmatrix}
        1-\beta_{uw} & \beta_u & 0 & \cdots & 0\\
        0 & 1-\beta_{uw} & \beta_{uw} & & 0\\
        \vdots & & \ddots & \ddots & 0\\
        0 & \cdots & \cdots & 1-\beta_{uw} & \beta_{uw}\\
        \beta_{uw} & 0 & \cdots & 0 & 1-\beta_{uw}
      \end{bmatrix},$$
\vspace{5mm}
$$c^u={}^t[\alpha_{uw},\;\; \alpha_{uw},\; \cdots,\; \alpha_{uw},\;\; \alpha_{uw}+n\beta_{uw}/d],$$
System~(\ref{eq-jeu}) is written:
\begin{equation}\label{gen-2}
  x_i^{k+1}=\min_{u\in\CU}\max_{w\in\CW}\{[M^{uw}x^k]_i+c^{uw}_i\},\;1\leq i\leq n.
\end{equation}
System~(\ref{gen-2}) is a dynamic programming
equation associated to a stochastic game problem. Let us denote by $h$ the
operator giving the dynamics~(\ref{gen-2}), i.e.
$h:\BR^n\to\BR^n$:
$$h_i(x^k)=\min_{u\in\CU}\max_{w\in\CW}\{[M^{uw}x^k]_i+c^{uw}_i\},\;1\leq i\leq n\;.$$
It is easy to see that $h$ is an additive 1-homogeneous and
monotone operator. We can also prove, by using similar arguments as
in Proposition~\ref{prop}, that the graph $\CG(h)$ is strongly
connected if and only if there exists $u$ in $\CU$ and $w\in\CW$
such that $\beta_{uw}\neq 0$ (i.e. $\beta_{uw}\in (0,1]$). Then,
taking this assumption, we apply Corollary~\ref{coro} and conclude
that the system~:
\begin{equation}\label{epd-app-3}
  \mu+x_i=\min_{u\in\CU}\max_{w\in\CW}\left\{(M^{uw}x)_i+c^{uw}_i\right\},\; 1\leq i \leq n.
\end{equation}
admits a solution $(\mu,x)$ where $x$ is defined up to an
additive constant, not necessarily in a unique way, and $\mu$ is
unique and satisfies:
\begin{equation}\nonumber
  \mu=\lim_{k\rightarrow +\infty}\frac{1}{k}\;x_i^k,\;\;\; 1\leq i \leq n\;.
\end{equation}
$\mu$ is interpreted as the average car speed.
\begin{theorem}\label{theo-gen-2}
  A solution $(\mu,x)$ of equation~(\ref{epd-app-3}) is given by~:
  $$\mu=\min_{u\in\CU}\max_{w\in\CW}\{\alpha_u+\frac{\beta_u}{d}\} \quad \text{and} \quad
  x= {}^t[0 \quad 1/d \quad 2/d \quad \cdots \quad (n-1)/d].$$
\end{theorem}
\proof Using the same arguments
as in the proof of Theorem~\ref{th-gen-2}, let
$(\bar{u},\bar{w})\in\CU\times\CW$ satisfying:
$$\mu=\min_{u\in\CU}\max_{w\in\CW}\{\alpha_{uw}+\frac{\beta_{uw}}{d}\}
    =\alpha_{\bar{u}\bar{w}}+\frac{\beta_{\bar{u}\bar{w}}}{d},$$
and let $s\in\CS$ be the feedback strategy given by
$(M^{\bar{u}\bar{w}},c^{\bar{u}\bar{w}})$, that is~:
$$s: x_i\rightarrow (\bar{u},\bar{w}),\; 1\leq i \leq n,$$
The couple $(\mu,x)$ given by~:
$$\mu=\alpha_{\bar{u}\bar{w}}+\frac{\beta_{\bar{u}\bar{w}}}{d}, \quad \text{and} \quad
  x= {}^t[0 \quad 1/d \quad 2/d \quad \cdots \quad (n-1)/d],$$
is solution of~:
  $$\mu+x=M^{\bar{u}\bar{w}}x+c^{\bar{u}\bar{w}}\; .$$
The strategy $s$ is optimal because for all $i$ in $\{1,2,\cdots,n\}$,
for all $u$ in $\CU$ and all $w$ in $\CW$ we have~:
$$[M^{\bar{u}\bar{w}}x+c^{\bar{u}\bar{w}}]_i = \mu+x_i
                    = \min_{u\in\CU}\max_{w\in\CW}\{\alpha_{uw}+\frac{\beta_{uw}}{d}\}+x_i
            = \min_{u\in\CU}\max_{w\in\CW}[M^{uw}x+c^{uw}]_i,$$
Hence the couple $(\mu,x)$ satisfies the spectral equation (\ref{epd-app-3}). \endproof

\begin{corollary}\label{cor3}
  The fundamental diagram on a circular road where the traffic is described by the dynamics
  (\ref{gen-2}) is given by~:
  $f=\min_{u\in\CU}\max_{w\in\CW}\{\alpha_{uw}d+\beta_{uw}\}$.
\end{corollary}
\proof Similarily, the average flow is equal to the average speed given in Theorem~\ref{theo-gen-2} mutiplied
by the average density $d$. \endproof

\section{Examples}

\subsection{Approximation of fundamental traffic diagrams}

\begin{figure}[htbp]
  \begin{center}
    \includegraphics[width=10cm]{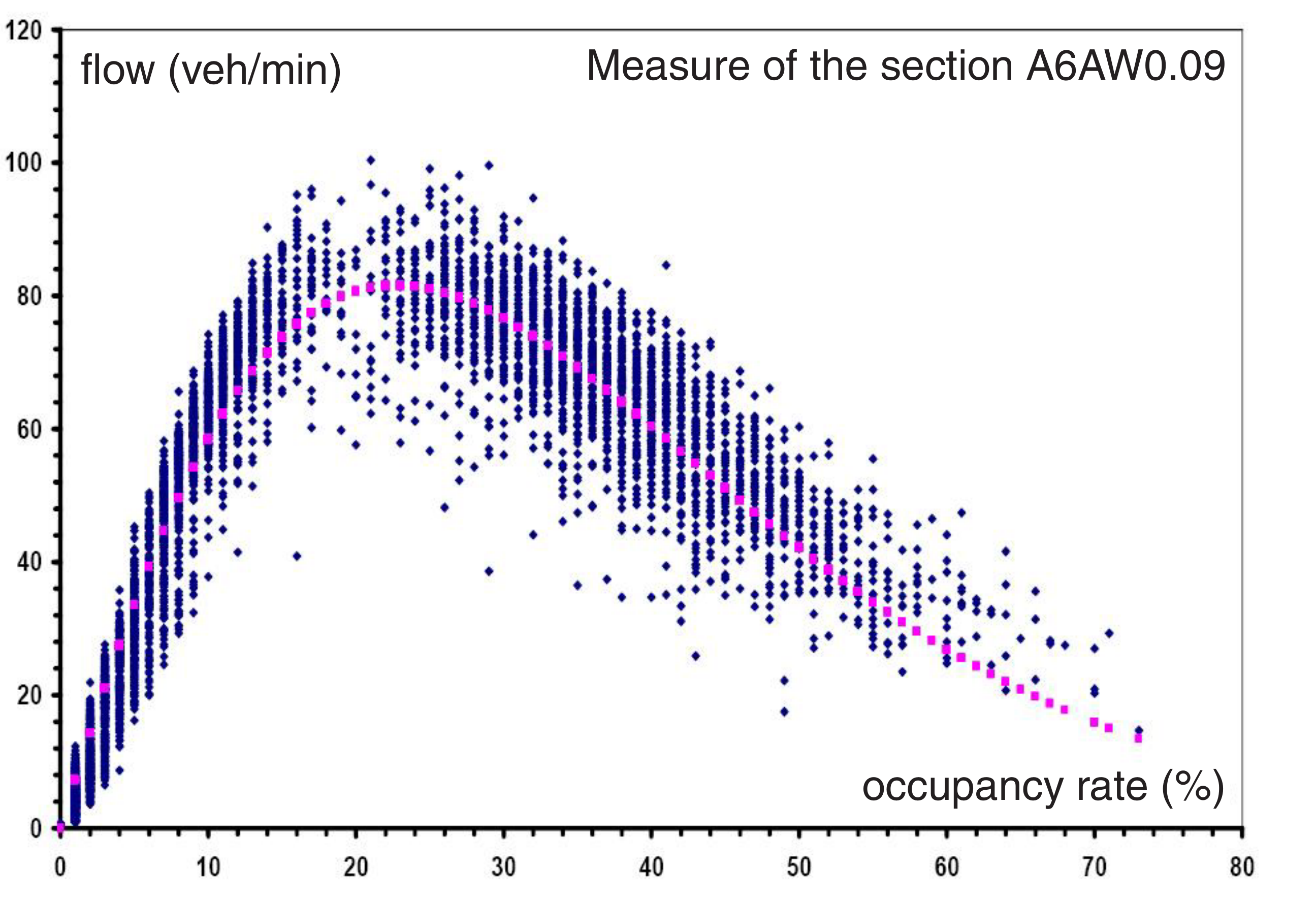}
    \caption{The fundamental diagram observed on the A6.}
    \label{approx0}
  \end{center}
\end{figure}

On Figure~\ref{approx0} we take an example of a fundamental
diagram obtained experimentally based on real measurements made on
a stretch of three lanes of the French highway~A6. We give an approximation
of this diagram using the stochastic game model.

On the x-axis of Figure~\ref{approx0} we have the car occupancy rate
on the road, which is a normalized density. On the
y-axis we have the car flow given by the number of
cars per minute. To obtain a normalized diagram where
the density and the flow of vehicles are given by relative
quantities in a coordinate system without unities, we normalized
the flow. To do this, we set arbitrarily to~1 the \emph{free
speed}, which is the average speed of vehicles in very low
densities. This quantity is given by the slope of the fundamental
diagram at the origin. Assuming that the maximum possible car flow
corresponds to the full density of vehicles moving freely (with
the free speed), we obtain the flow scale. For example if we
take $d=0.1$ as a very low density, then from Figure~\ref{approx0},
the flow corresponding to $d=10\%$ is $60$ veh./min., then we get
a maximum flow of $600$ veh./min (witch corrsponds to $d=100\%$). Then
by dividing the y-axis by $600$, we obtain a normalized diagram.

\begin{figure}[ht]
  \begin{center}
    \includegraphics[width=10cm,height=5cm]{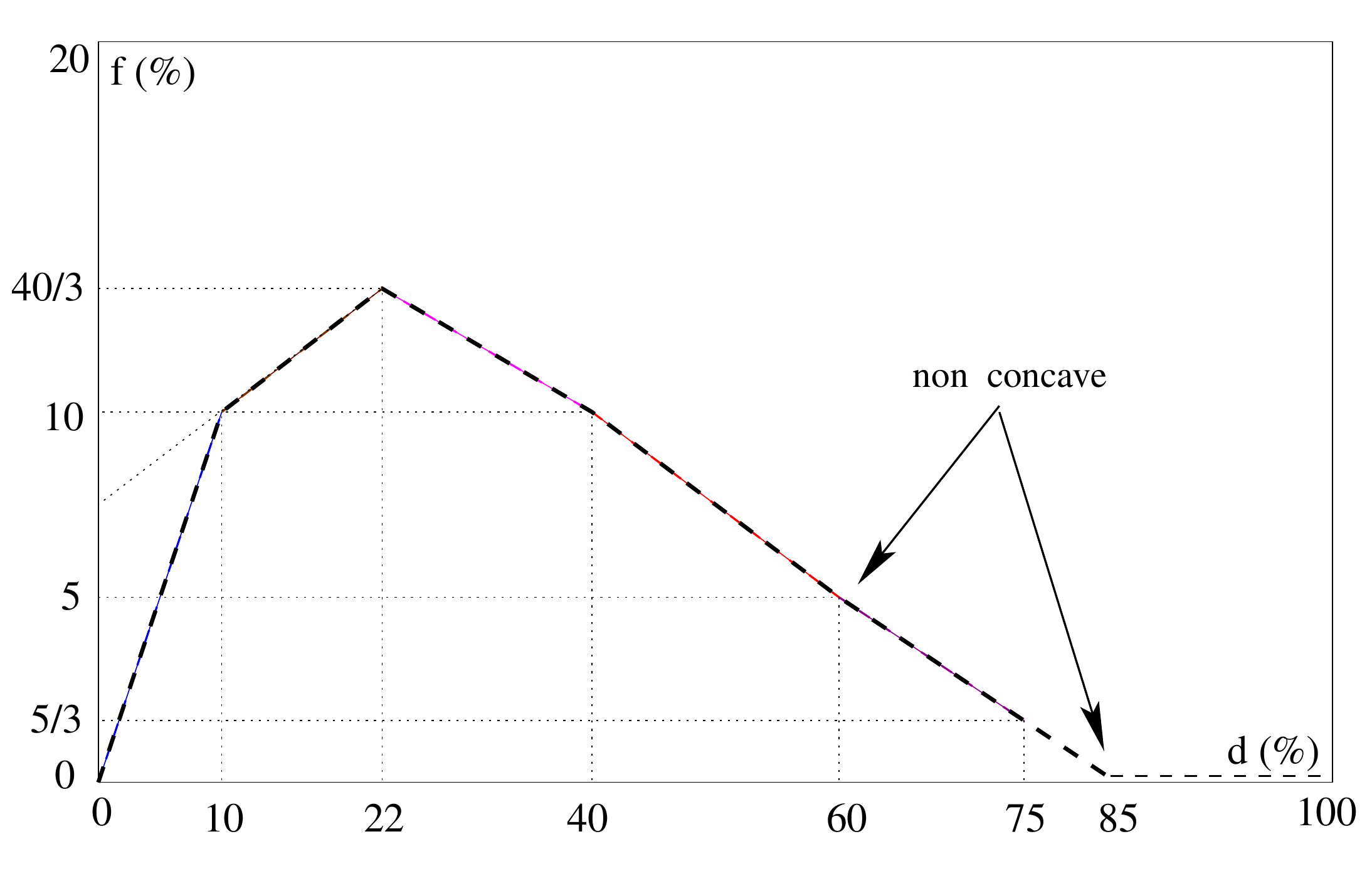}
    \caption{Stochastic game approximation.}
    \label{approxim3}
  \end{center}
\end{figure}

The objective here is to approximate the diagram of Figure~\ref{approx0},
in order interpret it and understand the traffic phases. To be able to approximate
non concave parts of the diagram, we use the stochastic game approximation.
Basing on Corollary~\ref{cor3}, we propose the following approximation (with six segments):
\begin{equation}\nonumber
  f=\min\left\{d, 0.27d+0.07, -0.19d+0.18, \max\{-0.25d+0.2,-0.2d+0.17,0\}\right\}\;,
\end{equation}
which is shown on Figure~\ref{approxim3}. 

\subsection{Traffic simulation and transitory regimes}

\begin{figure}[htbp]
    \begin{center}
      \includegraphics[width=2.5cm,height=2.5cm]{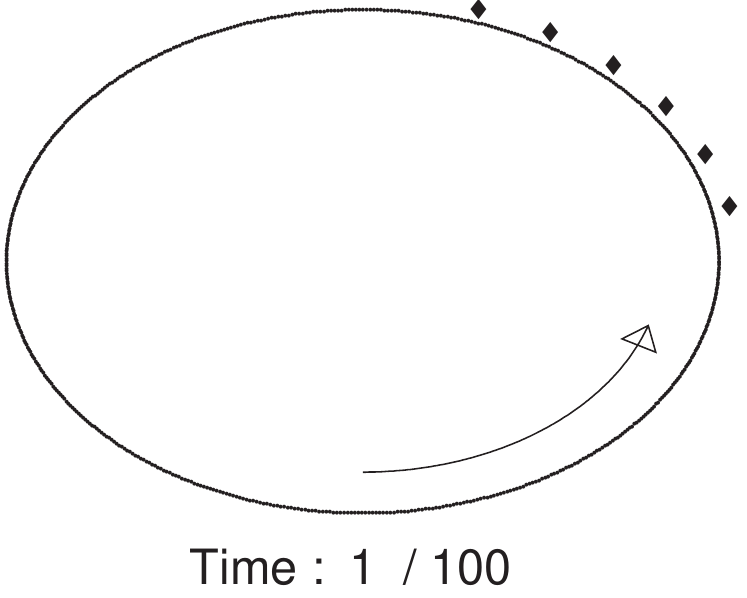}\hspace{8mm}
      \includegraphics[width=2.5cm,height=2.5cm]{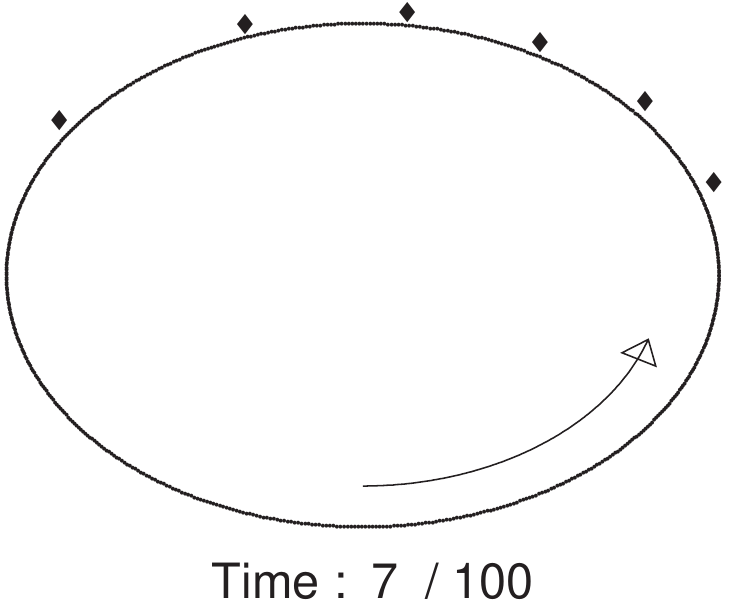}\hspace{8mm}
      \includegraphics[width=2.5cm,height=2.5cm]{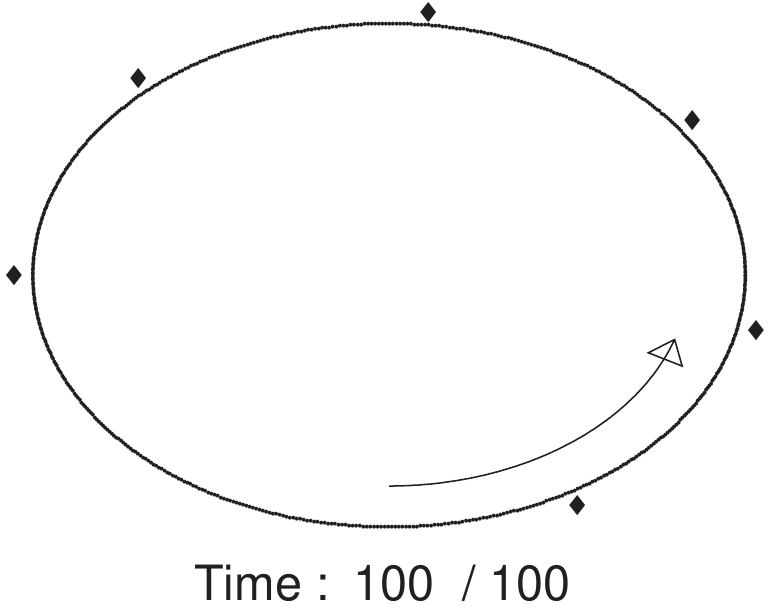}
      \caption{Phase 1.}
      \label{ph-1}
    \end{center}
    \label{DeuRou1}
\end{figure}

\begin{figure}[htbp]
    \begin{center}
      \includegraphics[width=2.5cm,height=2.5cm]{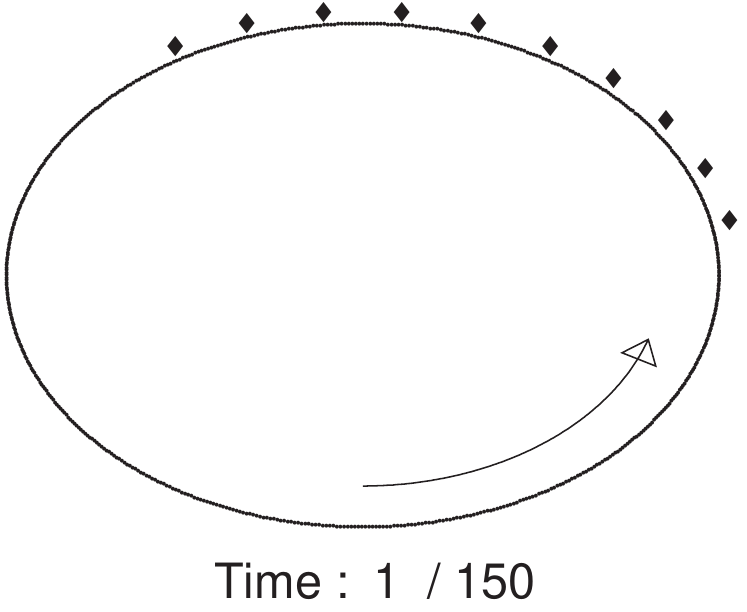}\hspace{8mm}
      \includegraphics[width=2.5cm,height=2.5cm]{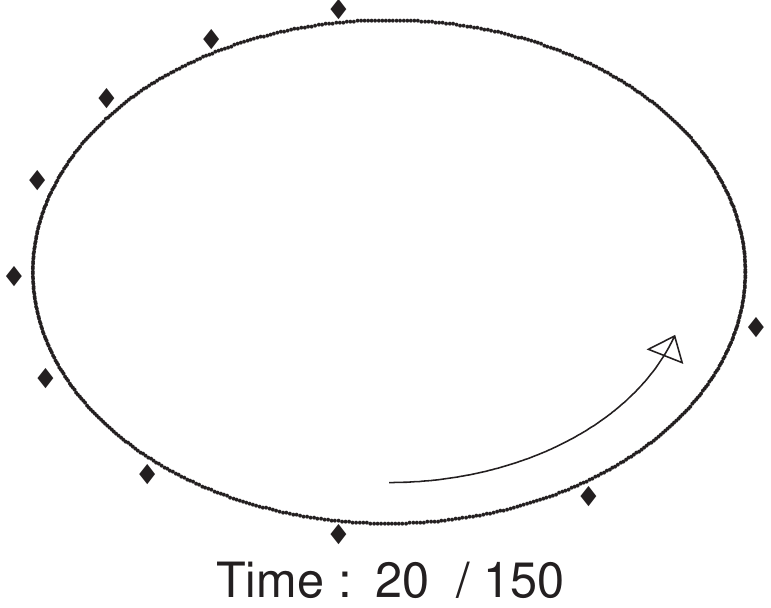}\hspace{8mm}
      \includegraphics[width=2.5cm,height=2.5cm]{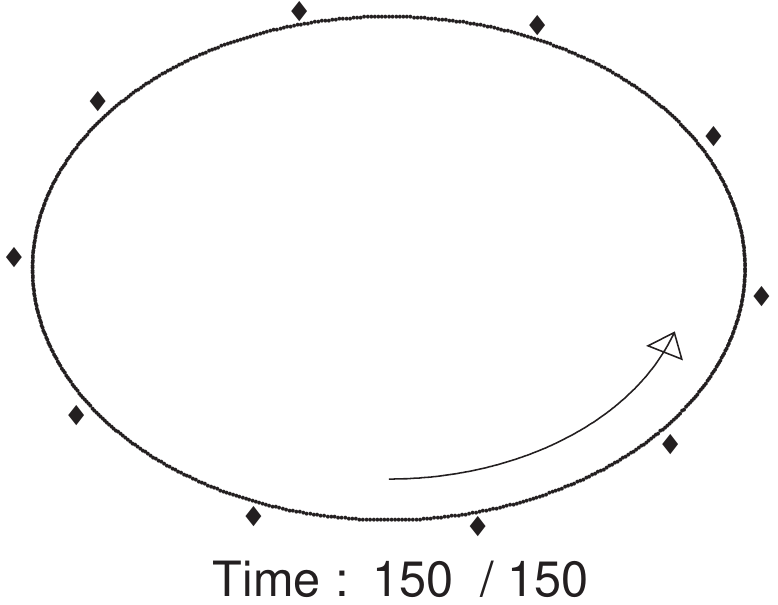}
      \caption{Phase 2.}
      \label{ph-2}
    \end{center}
    \label{DeuRou2}
\end{figure}

\begin{figure}[htbp]
    \begin{center}
      \includegraphics[width=2.5cm,height=2.5cm]{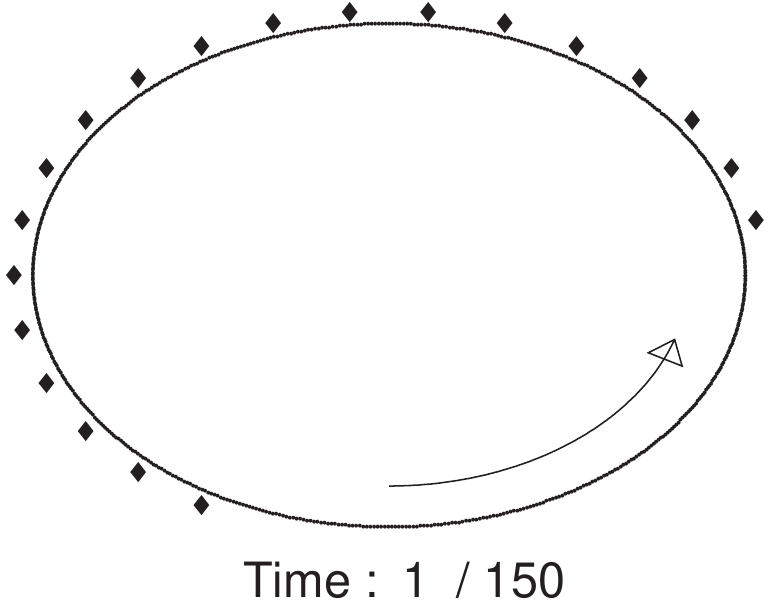}\hspace{8mm}
      \includegraphics[width=2.5cm,height=2.5cm]{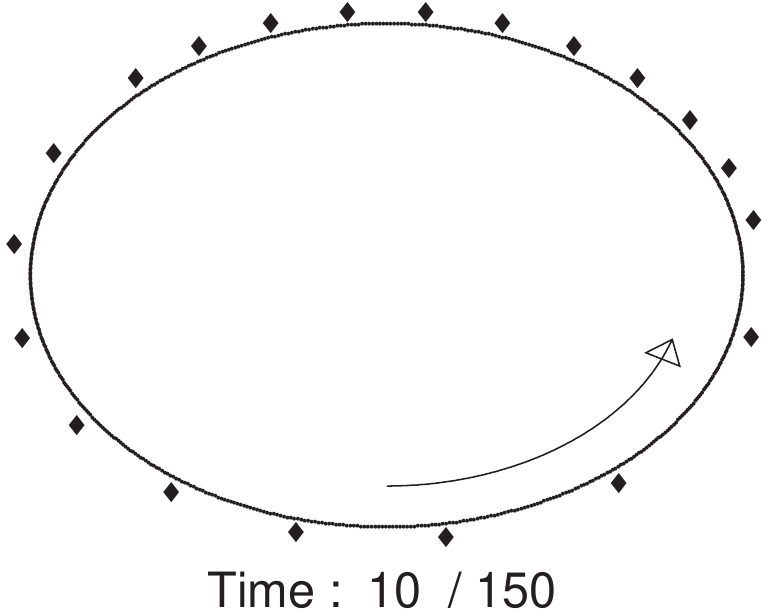}\hspace{8mm}
      \includegraphics[width=2.5cm,height=2.5cm]{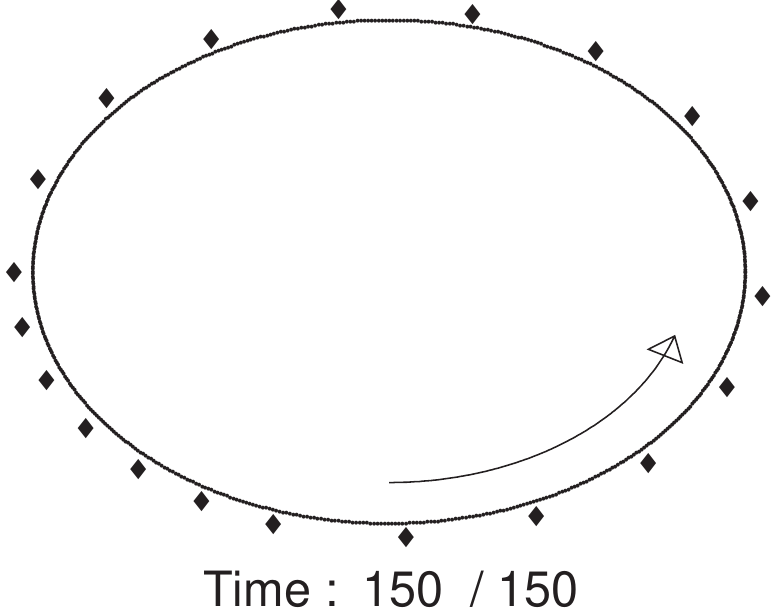}
      \caption{Phase 3.}
      \label{ph-3}
    \end{center}
    \label{DeuRou3}
\end{figure}

  Let us take $\CU=\{u_1,u_2,u_3\}$, with $(\alpha_1,\beta_1)=(1,0)$,
  $(\alpha_2,\beta_2)=(1/3,1/8)$, and $(\alpha_3,\beta_3)=(-1,1)$.
  The fundamental traffic diagram derived using the stochastic optimal control model is
  $f=\min\{d,\;(1/3)d +1/8,\;1-d\}$. The diagram presents three phases.
  On Figures~\ref{ph-1}, Figure~\ref{ph-2} and Figure~\ref{ph-3}, we simulate the
  traffic phases (phase 1: $f(d)=d$, phase 2: $f(d)=1/3 d+1/8$, and phase 3: $f(d)=1-d$).
  The car positions on the circular ring are given at three different times in
  order to highlight the transitory regime.

  We note that during phase~1 and phase~3, where $\beta_1=0$
  and $\beta_3=1$, the asymptotic car distributions on the road are not
  necessarily uniform, as shown on Figure~\ref{ph-1} and Figure~\ref{ph-3}.
  However, during phase~2 where $\beta_2=1/8\in (0,1)$, the asymptotic
  car distribution on the road is uniform, as obtained on Figure~\ref{ph-2}.

In general, and for the three models presented above, we observed numerically the following~:
\begin{itemize}
  \item For densities corresponding to the first segment of
      the fundamental diagram (the segment starting with the point $(0,0)$), which are
      in an interval of type $[0,d_0]$, the asymptotic car distribution
      on the road can be non uniform.
  \item If the last segment of the fundamental diagram (the segment ending with $(1,0)$)
      is given by $f(d)=1-d$ (the only case corresponding to
      $\beta_{\bar{u}}=1$), then for the corresponding densities,
      which are in an interval of type $[d_1,1]$, the asymptotic car distribution
      on the road can be non uniform.
  \item For all other densities $[d_0,d_1]$, the car distribution
      on the road converges to the uniform distribution.
\end{itemize}


\subsection*{General remarks}

\begin{itemize}
\item For the three models given in this article, the constraints
$\beta_u\in[0,1], \forall u\in\CU$ (or $\beta_{uw}\in[0,1],
\forall u\in\CU, \forall w\in\CW$) put the fundamental traffic
diagram in the triangle $[(0,0),(1,0),(0,1)]$ (triangle $ABD$ on
Figure~\ref{half}). Thus, the fundamental diagram lives on only
one half of the surface on witch it can a priori live, which is
the rectangle $[(0,0),(1,0),(1,1),(0,1)]$ (rectangle $ABCD$ on
Figure~\ref{half}).
\begin{figure}[htbp]
  \begin{center}
    \includegraphics[width=6cm,height=3.5cm]{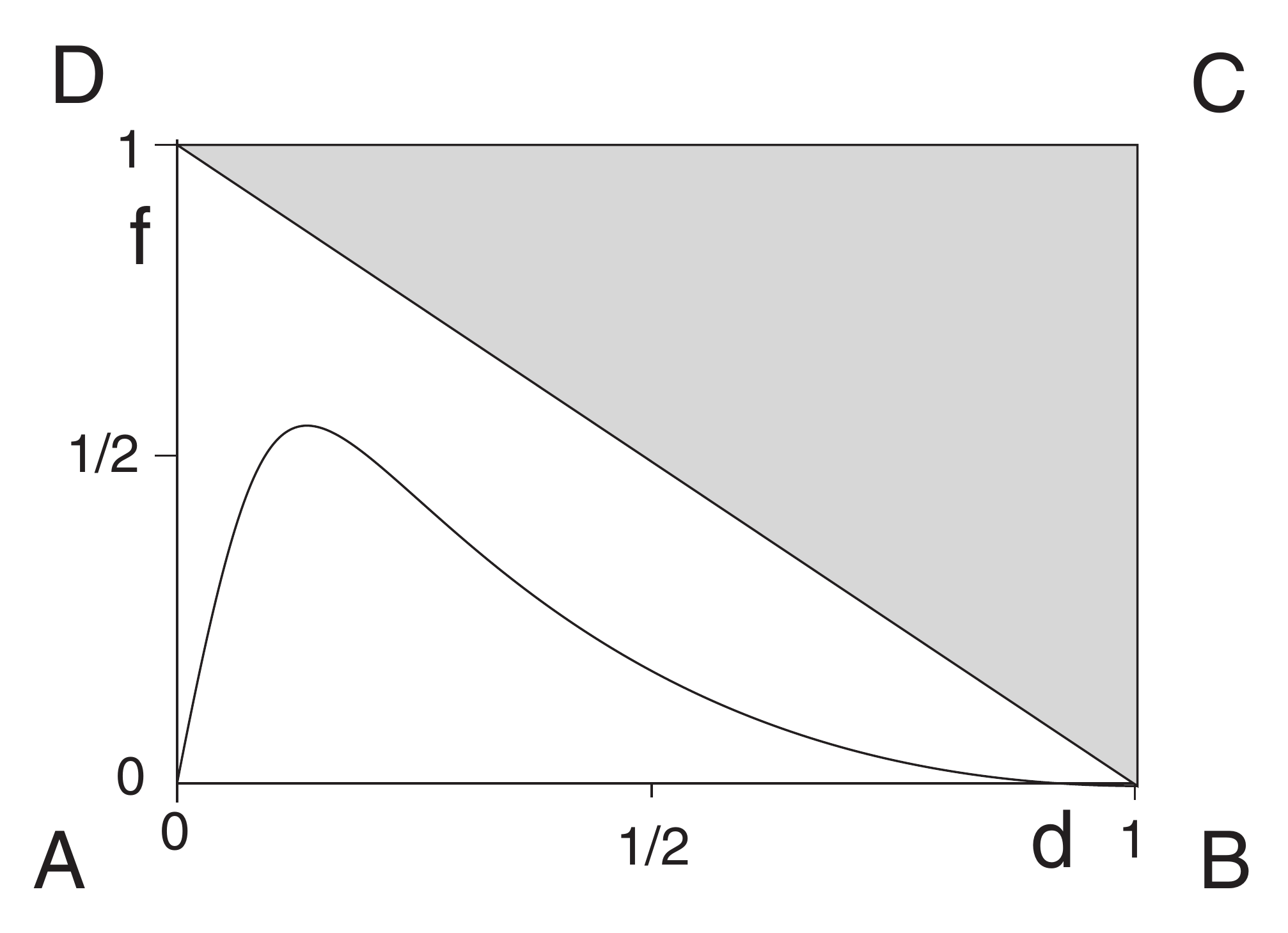}
    \caption{The fundamental diagram lives on the triangle $ABD$.}
    \label{half}
  \end{center}
\end{figure}
This can be written: $\forall d\in [0,1], f(d)\leq 1-d$. Indeed,
if we suppose that the fundamental diagram passes through a point
$p=(d,f(d))$ satisfying $f(d)>1-d$, then we can easily check that
with segments $\alpha_u d+\beta_u$ (or $\alpha_{uw}
d+\beta_{uw}$) satisfying $\beta_u\in[0,1]$ (or $\beta_{uw}\in[0,1]$),
we can never reach any point $p'=(d',f(d'))$ satisfying $d'\geq d$
and $f(d')\leq 1-d'$. Thus, the point $B=(1,0)$ is never reached.
This is absurd.

\item Finally, let us note that although the models given in this
article have stochastic interpretations, they are deterministic models.
A dual approach of the min-plus modeling given here is presented
in~\cite{LMQ05}, with an extension to a stochastic model.
The extension assumes that the desired speed
is stochastic and lives in a set of two reals. In this case, the
average speed is given by a Lyapunov exponent of a stochastic
min-plus matrix~\cite{BCOQ92,LMQ05}. In~\cite{Far08}, Petri net
and min-plus based models for the derivation of
fundamental diagrams of 2D-traffic are presented.
\end{itemize}

\section{Conclusion}

The models presented in this article describe 1D-traffic by dynamic
programming equations associated to
optimal control problems. By solving analytically these equations,
we derived explicitly the fundamental traffic diagrams. The derivation
permits the approximation of
a large class of fundamental diagrams. Moreover, the parameters
used in the models are basic traffic variables such as the desired velocity
of drivers, or the safety distance between
successive cars. This allows us to give simple interpretations
of the traffic phases appearing on experimental
fundamental diagrams obtained from real measurements.

\bibliographystyle{50}

\end{document}